\documentclass[3p]{elsarticle}

\usepackage{amssymb,amsmath,amsthm}

\usepackage[utf8]{inputenc}

\usepackage{tikz}
\usetikzlibrary{matrix,calc,arrows}
\tikzset{>=latex'}

\usepackage{hyperref}

\title{Spin structures on flat manifolds}

\author[rl]{Rafa{\l} Lutowski}
\ead{rafal.lutowski@mat.ug.edu.pl}

\author[bp]{Bartosz Putrycz}
\ead{bartosz.putrycz@mat.ug.edu.pl}

\address[rl,bp]{Institute of Mathematics, University of Gda\'nsk, ul. Wita Stwosza 57, 80-952 Gda\'nsk, Poland}

\newcommand{\GL}{\ensuremath{\mathrm{GL}}}
\newcommand{\SL}{\ensuremath{\mathrm{SL}}}
\newcommand{\G}{\ensuremath{\Gamma}}

\newcommand{\C}{\ensuremath{\mathbb C}}
\newcommand{\R}{\ensuremath{\mathbb R}}
\newcommand{\Q}{\ensuremath{\mathbb Q}}
\newcommand{\Z}{\ensuremath{\mathbb Z}}
\newcommand{\h}{\ensuremath{\mathbb H}}

\newcommand{\N}{\ensuremath{\mathbb N}}

\newcommand{\ses}[5]{#1 \longrightarrow #2 \longrightarrow #3 \longrightarrow #4 \longrightarrow #5}

\newcommand{\st}{\ensuremath{|\;}}

\newtheorem{thm}{Theorem}
\newtheorem{prop}[thm]{Proposition}
\newtheorem{lem}[thm]{Lemma}
\newtheorem{cor}{Corollary}
\theoremstyle{remark}
\newtheorem{rem}{Remark}
\theoremstyle{definition}
\newtheorem{df}{Definition}

\DeclareMathOperator{\vspan}{span}
\DeclareMathOperator{\Spin}{Spin}
\DeclareMathOperator{\SO}{SO}
\DeclareMathOperator{\GO}{O}
\DeclareMathOperator{\ind}{ind}
\DeclareMathOperator{\diag}{diag}

\NewDocumentCommand{\rl}{m}{\textbf{#1}}

\usepackage{subfig}

\begin{document}

\begin{abstract}
We present an algorithmic approach to the problem of the existence of spin structures on flat manifolds. We apply our method in the cases of flat manifolds of dimensions 5 and 6.
\end{abstract}


\maketitle

\section{Introduction}

Let $\G$ be an $n$ dimensional crystallographic group, i.e. a~discrete and cocompact subgroup of the group $E(n) = O(n) \ltimes \R^n$ of isometries of the Euclidean space $\R^n$. 
By the Bieberbach theorems (see \cite{Bi11, Bi12, Ch86}), $\G$ fits into short exact sequence
\begin{equation}
\label{eq:sesb}
\ses{0}{\Z^n}{\G}{G}{1},
\end{equation}
where $\Z^n$ is a maximal abelian normal subgroup of $\G$ and $G$ is a finite group, the so called holonomy group of $\G$. When in addition $\G$ is torsionfree, then $\G$ is called a Bieberbach group. In this case the orbit space $\R^n/\G$ is a flat manifold, i.e. a closed connected Riemannian manifold with sectional curvature equal to zero.

The existence of a spin structure on a manifold $X$ allows us to define on $X$ a Dirac operator. Every oriented flat manifold of dimension less than or equal to 3 admits a spin structure. In dimension 4, 24 out of 27 flat manifolds have spin structures (see \cite{PS10}). In this paper we present an algorithm to determine the existence of a spin structure on a flat manifold and present some facts concerning spin structures on flat manifolds of dimensions 5 and 6.

Section 2 recalls some basic definitions and introduces the necessary notations concerning Clifford algebras. The main goal of Section 3 is to present a more flexible form of a Pf\"affle criterion of the existence of spin structures on flat manifolds. The key tool in looking for spin structures on a flat manifold is the restriction of its holonomy representation to the Sylow $2$-subgroup of the holonomy group. In Section 4 we show that this restriction can be realized in a very convenient form and in Section 5 we show its usage in the criterion mentioned above. The algorithm for determining spin structures on flat manifolds is presented in Section 6 and is followed by an example of its usage for a 5-dimensional flat manifold. The last section presents some facts about spin structures for 5 and 6 dimensional manifolds.

\section{Clifford algebras and Spin groups}

\begin{df}
Let $n \in \N$. The \emph{Clifford algebra} $C_n$ is a~real associative algebra with one, generated by elements $e_1,\ldots,e_n$, which satisfy relations:
\[
\forall_{1 \leq i < j \leq n} \; e_i^2 = -1 \text{ and } e_ie_j = - e_je_i.
\]
\end{df}

\begin{rem}
We have the following $\R$-algebras isomorphisms:
\[
C_0 \cong \R, ~\ C_1 \cong \C, ~\ C_2 \cong \h.
\]
\end{rem}

\begin{rem}
We may view $\R^n := \vspan\{e_1, \ldots, e_n\}$ as a~vector subspace of $C_n$, for $n \in \N$.
\end{rem}

\begin{df}[Three involutions]
Let $n \in \N$. We have the following involutions of $C_n$:
\begin{itemize}
\item $^* \colon C_n \to C_n$, defined on the basis of (the vector space) $C_n$ by
\[ 
\forall_{1 \leq i_1 < i_2 < \ldots < i_k \leq n} \; ( e_{i_1}\ldots e_{i_k})^* = e_{i_k} \ldots e_{i_1};
\]
\item $' \colon C_n \to C_n$, defined on the generators of (the algebra) $C_n$ by 
\[
\forall_{1 \leq i \leq n} e_i' = -e_i.
\]
\item $\overline{ \phantom{a} } \colon C_n \to C_n$ -- the composition of the previous involutions
\[
\forall_{a \in C_n} \; \overline{a} = (a')^*.
\]
\end{itemize}
\end{df}

We are now ready to define the spin groups as subgroups of unit groups in the Clifford algebras:
\begin{df}
\[
\forall_{n \in \N} \;
\Spin(n) := \{ x \in C_n \st x' = x \wedge x\overline{x}=1 \}.
\]
\end{df}

\begin{prop}[{\cite[Prop. 6.1, page 86]{Sz12},\cite[page 16]{F00}}]
Let $n \in \N$. The map $\lambda_n \colon \Spin(n) \to \SO(n)$, defined by
\[
\forall_{x \in \Spin(n)} \forall_{v \in \R^n} \; \lambda_n(x)v = xv\overline{x}
\]
is a~continuous group epimorphism with kernel equal to $\{ \pm 1 \}$. Moreover for $n \geq 3$ $\Spin(n)$ is simply connected and $\lambda_n$ is the universal covering of $\SO(n)$.
\end{prop}

\section{Spin structures on (flat) manifolds}

\begin{df}
Let $X$ be an orientable closed manifold of dimension $n$. Let $Q$ be its principal $\SO(n)$-tangent bundle. A~\emph{spin structure} on $X$ is a~pair $(P,\Lambda)$, such that $P$ is a~principal $\Spin(n)$-bundle over $X$ and $\Lambda\colon P \to Q$ is a~$2$-fold covering for which the following diagram commutes:
\begin{center}
\begin{tikzpicture}
\node (x) at (2,0) {$X$};
\node (q) at (0,-1) {$Q$};
\node (p) at (0,1) {$P$};
\node (sp) at (-3,1) {$P \times \Spin(n)$};
\node (so) at (-3,-1) {$Q \times \SO(n)$};
\draw[->] (sp) -- (p);
\draw[->] (sp) --node[right]{$\Lambda \times \lambda_n$} (so);
\draw[->] (so) -- (q);
\draw[->] (p) --node[left]{$\Lambda$} (q);
\draw[->] (p) -- (x);
\draw[->] (q) -- (x);
\end{tikzpicture}
\end{center}
where the maps in the rows are defined by the action of the groups $\Spin(n)$ and $\SO(n)$ on the principal bundles $P$ and $Q$ respectively.
\end{df}

\begin{prop}[{\cite[page 40]{F00}}]
An orientable closed manifold $X$ has a spin structure if and only if its second Stiefel-Whitney class vanishes:
\[
w_2(X) = 0.
\]
Moreover in this case spin structures on $X$ are classified by $H^1(X,\Z_2)$.
\end{prop}

By the following proposition determining spin structures on flat manifolds becomes purely algebraic.

\begin{prop}[{\cite[Proposition 3.2]{P00} }]
Let $X$ be an $n$-dimensional orientable flat manifold with fundamental group $\G \subset E(n)$. Then the set of spin structures on $X$ is in bijection with the set of homomorphisms of the form $\varepsilon \colon \G \to \Spin(n)$ for which the following diagram commutes:
\begin{equation}
\label{eq:spin}
\begin{tikzpicture}[baseline=(current  bounding  box.center)]
\node (g) at (0,0) {$\G$};
\node (sp) at (2,1.5) {$\Spin(n)$};
\node (so) at (2,0) {$\SO(n)$}; 
\draw[->] (g) --node[above left]{$\varepsilon$} (sp);
\draw[->] (sp) -- node[right]{$\lambda_n$} (so);
\draw[->] (g) -- node[above]{$r$} (so);
\end{tikzpicture}
\end{equation}
where $r \colon \G \to \SO(n)$ takes the rotational part of every element of $\G$:
\[
\forall_{(A,a) \in \G \subset \SO(n) \ltimes \R^n} r(A,a) = A.
\]
\end{prop}

\begin{rem}
By a little abuse of notation we will call $r \colon \G \to \SO(n)$ the \emph{holonomy representation} of $\G$.
\end{rem}

Now		 let $X=\R^n/\G$ be an orientable flat manifold with fundamental group $\G \subset E(n)$. The group $\G$ is finitely presented. Let
\[
\G = \langle S \;|\; R \rangle
\]
be its presentation with the set of generators $S$ closed under taking inversions ($S^{-1}=S$) and the set of relations $R$, both finite sets. A map $\varepsilon\colon S \to \Spin(n)$ can be extended to a homomorphism $\varepsilon \colon \G \to \Spin(n)$ if and only if it preserves the relations of $\G$:
\[
\forall_{r_1,\ldots,r_l \in S} r_1 \cdot \ldots \cdot r_l \in R \Rightarrow \varepsilon(r_1) \cdot \ldots \cdot \varepsilon(r_l) = 1.
\]
Moreover, since $\ker \lambda_n = \{\pm 1\}$, in order to get commutativity of the diagram \eqref{eq:spin} we must have
\[
\forall_{\gamma \in \G} \forall_{x \in \Spin(n)} r(\gamma) = \lambda_n(x) \Rightarrow \varepsilon(\gamma) = x \vee \varepsilon(\gamma) = -x.
\]

Hence to check if we can construct a homomorphism $\varepsilon$ then for every generator $s \in S$ it is enough to find an element $x \in \Spin(n)$ such that
\[
r(s) = \lambda_n(x)
\]
and check which combinations of signs of those elements of $\Spin(n)$ preserve relations of $\G$.

In general it is not an easy task to find preimages $\lambda_n^{-1}(g)$ of an element $g \in \SO(n)$. The following proposition allows us to search for such finite subgroups of $\SO(n)$ which are easier to work with.

\begin{prop}[{\cite[Proposition 2.1]{HiSz08}}]
Let $n \in \N$ and $\G_1, \G_2 \subset E(n)$ be isomorphic Bieberbach groups. Then the set of spin structures on the flat manifold $\R^n/\G_1$ is in bijection with the set of spin structures on the flat manifold $\R^n/\G_2$.
\end{prop}

The bijection in the above proposition is given as follows. Let $\G_1,\G_2$ be Bieberbach groups as above. Let $r_i\colon \G_i \to SO(n)$ be the holonomy representations and let $G_i=r_i(\G_i)$ be the holonomy groups of $\G_i$, for $i=1,2$. If $\varepsilon_1 \colon \G_1 \to \Spin(n)$ defines a spin structure on $\R^n/\G_1$ then the corresponding homomorphism $\varepsilon_2 \colon \G_2 \to \Spin(n)$ fits into the following commutative diagram
\begin{equation}
\begin{tikzpicture}[baseline=(current  bounding  box.center)]
\node (g1) at (0,0) {$\G_1$};
\node (p1) at (2,1) {$G_1$};
\node (s1) at (2,3) {$\varepsilon_1(\G_1)$};
\node (g2) at (4,0) {$\G_2$};
\node (p2) at (6,1) {$\SO(n)$};
\node (s2) at (6,3) {$\Spin(n)$};
\draw[->] (p1) -- node[above]{$\varphi$} (p2.west);
\draw[->] (g1) -- node[above]{$\Phi$} (g2);
\draw[->] (s1) -- node[above]{$\alpha$} (s2);
\draw[->] (g1) -- node[above left]{$\varepsilon_1$} (s1);
\draw[->] (g1) -- node[above]{$r_1$} (p1);
\draw[->] (s1) -- node[right]{$\lambda_n$} (p1);
\draw[line width=.6em,white] (g2) -- (s2);
\draw[->] (g2) -- node[above left]{$\varepsilon_2$} (s2);
\draw[->] (g2) -- node[above]{$r_2$} (p2);
\draw[->] (s2) -- node[right]{$\lambda_n$} (p2);
\end{tikzpicture}
\end{equation}
where $\Phi \colon \G_1 \to \G_2$ is the isomorphism, $\varphi$ is the homomorphism induced by $\Phi$. The map $\alpha$ is induced by $\varphi$ as follows: if $\varphi$ is a~conjugation by a~matrix $F \in \GL(n,\R)$ then $\alpha$ is a conjugation by a~lift $\tilde{F}$ of $F$ in the metalinear group $\operatorname{ML}(n,\R)$ -- a double cover of $\GL(n,\R)$.

\begin{cor}
\label{cor:spin}
Let $\G \subset E(n)$ be a Bieberbach group with holonomy representation $r \colon \G \to \SO(n)$ and holonomy group $G=r(\G)$. The set of spin structures on the flat manifold $\R^n/\G$ is in bijection with the set of homomorphisms of the form $\varepsilon\colon \G \to \Spin(n)$ for which the following diagram commutes:
\begin{equation}
\label{eq:spin1}
\begin{tikzpicture}[baseline=(current  bounding  box.center)]
\node (g) at (0,1.5) {$\G$};
\node (sp) at (2,1.5) {$\Spin(n)$};
\node (p) at (0,0) {$G$};
\node (so) at (2,0) {$\SO(n)$}; 
\draw[->] (g) --node[above]{$\varepsilon$} (sp);
\draw[->] (g) --node[right]{$r$} (p);
\draw[->] (sp) -- node[right]{$\lambda_n$} (so);
\draw[->] (p) -- node[above]{$\varphi$} (so);
\end{tikzpicture}
\end{equation}
where $\varphi \colon G \to \SO(n)$ is a representation of $G$ equivalent to the identity map $id \colon G \to G \subset \SO(n)$.
\end{cor}

\section{Flat manifolds with 2-group holonomy}

\begin{prop}[{\cite[Proposition 1]{DSS04} }]
Let $n \in \N$. Let $\G \in E(n)$ be a Bieberbach group with holonomy representation $r \colon \G \to \SO(n)$ and holonomy group $G=r(\G)$. Let $S \subset G$ be a $2$-Sylow subgroup of $G$. Then the flat manifold $\R^n/\G$ admits a spin structure if and only if $\R^n/r^{-1}(S)$ admits one.
\end{prop}

By Corollary \ref{cor:spin} in the process of determining the existence of spin structures on a flat manifold we can choose any subgroup of $\SO(n)$ which is conjugated in $\GL(n,\R)$ to its holonomy group. By the above proposition it is enough to look on $2$-subgroups of $\SO(n)$. In this section we will show that for every $2$-group in $\SO(n)$ we can find its conjugate $G$ in such a way that $\lambda_n^{-1}(G)$ is easy to compute.

\begin{rem}
The extension \eqref{eq:sesb} defines the \emph{integral holonomy representation} $\rho \colon G \to \GL(n,\Z)$ defined by the conjugations in $\G$:
\[
\forall_{g \in G} \forall_{z \in \Z^n} \rho_g(z) = \gamma z \gamma^{-1},
\]
where $\gamma$ is an element of $\G$ such that $r(\gamma) = g$. This representation is $\R$-equivalent to the "identity representation" $id \colon G \to G \subset \SO(n)$.
\end{rem}

\begin{thm}[{\cite[Theorem 1.10]{EM79}}]
Let $G$ be a finite $p$-group and let $\varphi\colon G \to \GL(m,\Q)$ be an irreducible representation over $\Q$. Then either $\varphi$ is induced from a representation of a subgroup of index $p$ or $[G:\ker \varphi] \leq p$.
\end{thm}

By an induction argument we immediately get

\begin{cor}
\label{cor:rep2grp}
Every irreducible rational representation of $2$-group is induced from a rational representation of degree $1$.
\end{cor}

Now let's take a closer look on a matrix representation of a $2$-group $G$
\[
\varphi \colon G \to \GL(n,\Q).
\]
By Corollary \ref{cor:rep2grp} we may assume that
\[
\varphi = \ind_{H_1}^G \tau_1 \oplus \ldots \oplus \ind_{H_k}^G \tau_k
\]
where $H_i$ is a subgroup of $G$ and $\tau_i \colon H \to \Q^*$ is a representation of $H_i$ of degree $1$, for $i=1,\ldots,k$. Since for every $1 \leq i \leq k$ we have $\tau_i(H_i)\subset\{ \pm 1 \}$, hence every element of $\varphi(G)$ is an orthogonal integral matrix and $\varphi$ is if the form
\[
\varphi \colon G \to \GO(n,\Z) := \GO(n) \cap \GL(n,\Z).
\]

Now if a $2$-group $G \subset \SO(n)$ is a holonomy group of a Bieberbach group $\G \subset E(n)$ then by Corollary \ref{cor:spin} the set of spin structures of the manifold $\R^n/\G$ is in bijection with the set of homomorphisms $\varepsilon \colon \G \to \Spin(n)$ which make the following diagram commute

\begin{equation}
\label{eq:spin2}
\begin{tikzpicture}[baseline=(current  bounding  box.center)]
\node (g) at (0,1.5) {$\G$};
\node (sp) at (2,1.5) {$\Spin(n)$};
\node (p) at (0,0) {$G$};
\node (so) at (2,0) {$\SO(n,\Z)$}; 
\draw[->] (g) --node[above]{$\varepsilon$} (sp);
\draw[->] (g) --node[right]{$r$} (p);
\draw[->] (sp) -- node[right]{$\lambda_n$} (so);
\draw[->] (p) -- node[above]{$\varphi$} (so);
\end{tikzpicture}
\end{equation}
where $\SO(n,\Z) = \SL(n,\Z) \cap \SO(n)$. This seems to be a minor change in comparison to Corollary \ref{cor:spin}, but it simplifies a lot the problem of determining preimages of $\lambda_n$.

\section{Special orthogonal group over the integers}

In this section we will show how to determine the preimage of any element of the group $\SO(n,\Z)$ under the homomorphism $\lambda_n$ for $n\geq 2$. Recall that in this case $\ker \lambda_n = \{ \pm 1 \}$ so calculation of one element in the preimage immediately gives us the other one. 

The group $\GO(n,\Z)$ fits into the following exact sequence
\[
\ses{1}{N}{\GO(n,\Z)}{S_n}{1},
\] 
where $S_n$ is the symmetric group on $n$ letters and $N \subset \GO(n,\Z)$ is the group of diagonal matrices with $\pm 1$ on the diagonal. The sequence splits and the splitting homomorphism sends a permutation $\sigma \in S_n$ to its permutation matrix $P_\sigma \in \GO(n,\Z)$.

Now let $X \in \SO(n,\Z)$ be an integral orthogonal matrix. There exist inversions $\sigma_1, \ldots \sigma_k \in S_n$ and a diagonal integral matrix $D \in N$ such that
\[
X = P_{\sigma_1} \cdots P_{\sigma_k} \cdot D.
\]
Unfortunately matrices of inversions have determinant equal to $-1$ and they don't belong to $\SO(n,\Z)$. A little modification changes this fact. Let $(p\;q) \in S_n$ be an inversion with $p < q$. Define the matrix $P'_{(p\;q)} \in \SO(n,\Z)$ as follows:
\[
P'_{(p\;q)} = \diag(\underbrace{1,\ldots,1}_{p-1},-1,1,\ldots,1) \cdot P_{(p\;q)}.
\]
We get that
\begin{equation}
\label{eq:decomposition}
X =  P'_{\sigma_1} \cdots P'_{\sigma_k} \cdot D'
\end{equation}
where $D' \in N$ but this time all the factors $P'_{\sigma_1}, \ldots, P'_{\sigma_k},  D'$ in the decomposition of $X$ have determinant $1$ and hence they are elements of $\SO(n,\Z)$. In order to determine $\lambda_n^{-1}(X)$ it is enough to calculate the preimages of its factors:
\begin{lem} \ 
\begin{enumerate}
\item If $D' \in N$ is a matrix with $-1$ on the diagonal entries $n_1,\ldots,n_l$ ($l$ even) then
\begin{equation}
\label{eq:map_diag}
\lambda_n(\pm e_{n_1} \cdots e_{n_l}) = D'.
\end{equation}
\item If $(p\;q) \in S_n$ is an inversion with $p < q$ then
\begin{equation}
\label{eq:map_inv}
\lambda_n\left( \pm\frac{1+e_p e_q}{\sqrt{2}} \right) = P'_{(p\;q)}.
\end{equation}
\end{enumerate}
\end{lem}

\begin{proof} \ 
\begin{enumerate}
\item Let $i \in \{ n_1, \ldots, n_l \}$. Changing the order of the factors of the product $e_{n_1} \cdots e_{n_l}$ does not affect the value of $\lambda_n$. Hence without lose of generality we can assume that $i=n_1$. We get
\[
\begin{split}
\lambda_n(\pm e_{n_1} \cdots e_{n_l}) e_{n_1} & = \lambda_n(e_{n_1} \cdots e_{n_l}) e_{n_1} = e_{n_1} \cdots e_{n_l} e_{n_1} e_{n_l} \cdots e_{n_1}\\
& = -e_{n_1} \cdots e_{n_{l-1}} e_{n_1} e_{n_l} e_{n_l} \cdots e_{n_1} = \ldots = (-1)^{l-1} e_{n_1} e_{n_1} e_{n_2} \cdots e_{n_l} e_{n_l} \cdots e_{n_1}\\
& = (-1)^{l} e_{n_2} \cdots e_{n_{l-1}} (-1) e_{n_{l-1}} \cdots e_{n_2} e_{n_1} = \ldots =  (-1)^{2l-1} e_{n_1} = -e_{n_1} = D' e_{n_1}.
\end{split}
\]
On the other hand if $i \not\in \{n_1, \ldots, n_l \}$ then
\[
\begin{split}
\lambda_n(\pm e_{n_1} \cdots e_{n_l}) e_i & = e_{n_1} \cdots e_{n_l} e_{i} e_{n_l} \cdots e_{n_1} \\
& = -e_{n_1} \cdots e_{n_{l-1}} e_{i} e_{n_l} e_{n_l} \cdots e_{n_1} = \ldots = (-1)^l e_i e_{n_1} \cdots e_{n_l} e_{n_l} \cdots e_{n_1}\\
& = (-1)^{2l} e_i = e_i = D'e_i.
\end{split}
\]
\item We get
\[
\begin{split}
\lambda_n\left( \pm\frac{1+e_p e_q}{\sqrt{2}} \right) e_p & = \frac{1+e_p e_q}{\sqrt{2}} \cdot e_p \cdot  \frac{1+e_q e_p}{\sqrt{2}} = \frac{(e_p+e_pe_qe_p)(1+e_qe_p)}{2}\\
& = \frac{(e_p+e_q)(1+e_qe_p)}{2} = \frac{e_p+e_q+e_pe_qe_p+e_q^2e_p}{2} = e_q = P'_{(p\;q)} e_p.
\end{split}
\]
Similarly
\[
%\begin{split}
%\lambda_n\left( \pm\frac{1+e_p e_q}{\sqrt{2}} \right) e_q & = \frac{1+e_p e_q}{\sqrt{2}} \cdot e_q \cdot \frac{1+e_q e_p}{\sqrt{2}} = \frac{(e_q+e_pe_q^2)(1+e_qe_p)}{2}\\
%& = \frac{(e_q-e_p)(1+e_qe_p)}{2} = \frac{e_q-e_p+e_q^2e_p-e_pe_qe_p}{2} = -e_p = P'_{(p\;q)} e_q.
%\end{split}
\lambda_n\left( \pm\frac{1+e_p e_q}{\sqrt{2}} \right) e_q = -e_p = P'_{(p\;q)} e_q.
\]
If $i \not\in \{ p,q \}$ then the elements $(1+e_pe_q)/\sqrt{2}$ and $e_i$ commute, hence
\[
\lambda_n\left( \pm\frac{1+e_p e_q}{\sqrt{2}} \right) e_i = \frac{1+e_p e_q}{\sqrt{2}} \cdot e_i \cdot \frac{1+e_q e_p}{\sqrt{2}} = e_i \frac{1+e_p e_q}{\sqrt{2}} \cdot \frac{1+e_q e_p}{\sqrt{2}} = e_i = P'_{(p\;q)} e_i.
\]
\end{enumerate}
\end{proof}

\section{Notes about the algorithm}

Let $n \geq 2$. Assume that $\G' \subset E(n)$ is a Bieberbach group with holonomy representation $r \colon\G' \to \SO(n)$ and that $\G'$ fits into the following short exact sequence
\[
0 \longrightarrow \Z^n \stackrel{i}{\longrightarrow} \G' \stackrel{r}{\longrightarrow} G' \longrightarrow 1.
\]
The following steps will determine the existence of spin structures on the flat manifold $\R^n/\G'$.

\begin{description}
\item[Step 1] Determine a Sylow $2$-subgroup $G$ of $G'$ and its preimage $\G = r^{-1}(G)$ in $\G'$. We get an extension
\[
0 \longrightarrow \Z^n \stackrel{}{\longrightarrow} \G \stackrel{r}{\longrightarrow} G \longrightarrow 1
\]
where $r$ is in fact a restriction $r_{|\G}$.

\item[Step 2] Determine a representation $\varphi \colon G \to \SO(n,\Z)$ of a 2-group $G$ which is $\R$-equivalent to the identity representation $id_G \colon G \to G \subset \SO(n)$. Note that it may be helpful to build a list of all $\Q$-irreducible integral and orthogonal representations of $G$. Since we are in characteristic zero the character theory is very useful in determining which of those are subrepresentations of $id_G$.

\item[Step 3] Fix a generating set $\{ g_1, \ldots, g_s \}$ of $G$. For every $1 \leq i \leq s$ decompose $\varphi(g_i)$ as in \eqref{eq:decomposition} and then, using formulas \eqref{eq:map_diag} and \eqref{eq:map_inv} determine $x_i \in \Spin(n)$ such that 
\[
\lambda_n(x_i) = \varphi(g_i).
\]

\item[Step 4] Determine the integral holonomy representation
\[
\varrho \colon G \to \GL(n,\Z).
\]
Denote by $\varrho_{i,j}(g) \in \Z$ the entry in the $i$-th row and $j$-th column of the matrix $\varrho(g) \in \GL(n,\Z)$ where $1 \leq i,j \leq n, g \in G$. It is worth to notice that CARAT uses the integral holonomy representation to store crystallographic groups as a subgroup of $\GL(n,\Z) \ltimes \Q^n$ with translation lattice being always $\Z^n$. In this form the projection on the first coordinate defines the integral holonomy representation. Note that this is not a constraint in any way, since $\varrho, id_G$ and $\varphi$ are all $\R$-equivalent.

\item[Step 5] Let $a_1,\ldots,a_n \in \G$ be the images of the generators of $\Z^n$ in $\G$. Let $\gamma_1, \ldots, \gamma_s$ be elements of $\G$ such that
\[
\forall_{1 \leq i \leq s} r(\gamma_i) = g_i.
\]
By \cite[Proposition 1, page 139]{J97} 
\[
\G = \langle a_1, \ldots, a_n, \gamma_1, \ldots, \gamma_s \rangle.
\]
Note that if we have a homomorphism $\varepsilon \colon \G \to \Spin(n)$ such that $\lambda_n \varepsilon = r$ then
\[
\varepsilon(a_i) \in \{ \pm 1 \} \wedge \varepsilon(\gamma_j) \in \{ \pm x_i \}
\]
for all $1 \leq i \leq n, 1 \leq j \leq s$. Now for every possible value of a function $\varepsilon$ on the generators of $\G$ we have to check whether we can extend it to a homomorphism of groups, i.e. we have to check whether the images preserve the relations amongst the generators of $\G$ which are of three types:
\begin{enumerate}
\item Relations which come from the monomorphism $\Z^n \to \G$ are the commutator relations and they are automatically satisfied, since all the generators $a_1, \ldots, a_n$ are mapped to the center of $\Spin(n)$.
\item Relations which come from the action of $G$ on $\Z^n$. Let $1 \leq i \leq n$ and $1 \leq j \leq s$. Using the holonomy representation $\varrho$ we get the following relation in $\G$:
\[
\gamma_j a_i \gamma_j^{-1} = a_1^{\varrho_{1i}(g_j)} \cdots a_n^{\varrho_{ni}(g_j)}.
\]
The corresponding relation in $\Spin(n)$ should be as follows
\[
\varepsilon(a_1)^{\varrho_{1i}(g_j)} \cdots \varepsilon(a_n)^{\varrho_{ni}(g_j)} = \varepsilon(\gamma_j) \varepsilon(a_i) \varepsilon(\gamma_j)^{-1} = \varepsilon(a_i)
\]
since $\varepsilon(a_i) = \pm 1$. From the same reason the above equation may be written as
\[
\varepsilon(a_1)^{\varrho_{1i}(g_j)} \cdots \varepsilon(a_n)^{\varrho_{ni}(g_j)} \varepsilon(a_i) = 1.
\]
\item Relations which come from relations of $G$. Let
\[
g_{i_1} \cdots g_{i_k}
\]
be a relator of $G$ (you can skip inverses since $G$ is finite). Then
\[
\gamma_{i_1} \cdots \gamma_{i_k} = a_1^{\alpha_1} \cdots a_n^{\alpha_n}
\]
for some $\alpha_1, \ldots, \alpha_n \in \Z$. The resulting relation in $\Spin(n)$ is
\[
\varepsilon(\gamma_{i_1}) \cdots \varepsilon(\gamma_{i_k}) = \varepsilon(a_1)^{\alpha_1} \cdots \varepsilon(a_n)^{\alpha_n},
\]
which is equivalent to
\[
\varepsilon(\gamma_{i_1}) \cdots \varepsilon(\gamma_{i_k}) \varepsilon(a_1)^{\alpha_1} \cdots \varepsilon(a_n)^{\alpha_n} = 1.
\]
\end{enumerate}

\end{description}

\section{Example}

Let $\G'$ be a Bieberbach group generated by the matrices
\[
\begin{bmatrix}
 1 &  0 &   0 &  0 &  0 &  -1/3\\
 0 &  0 &  -1 &  0 &  0 &     0\\
 0 &  1 &  -1 &  0 &  0 &     0\\
 0 &  0 &  -1 &  1 &  0 &     0\\
 0 &  0 &  -1 &  0 &  1 &     0\\
 0 &  0 &   0 &  0 &  0 &     1\\
\end{bmatrix},
\begin{bmatrix}
 -1 &   0 &   0 &  0 &  0 &     0\\
  0 &   0 &  -1 &  1 &  1 &   1/2\\
  0 &  -1 &   0 &  1 &  1 &     0\\
  0 &   0 &   0 &  1 &  0 &  -1/2\\
  0 &   0 &   0 &  0 &  1 &     0\\
  0 &   0 &   0 &  0 &  0 &     1\\
\end{bmatrix}
\]
and the matrices of the form
\begin{equation}
\label{eq:lattice}
a_i = 
\begin{bmatrix}
I & e_i\\
0 & 1 \\
\end{bmatrix}
\end{equation}
where $I$ is the identity matrix of degree $5$ and the vectors $e_i, i=1,\ldots,5$ are generators of $\Z^5$. The group is denoted in CARAT by min.134.1.2.2. The holonomy group $G'$ of $\G'$ is isomorphic to the symmetric group $S_4$ and so its $2$-Sylow subgroup $G$ is isomorphic to the dihedral group $D_8$. If $r \colon \G' \to \SO(5)$ is the holonomy representation, then the preimage $\G = r^{-1}(G)$ is generated by the matrices $a_1, \ldots, a_5$ and the following ones:
\[
A = 
\begin{bmatrix}
 1 &   0 &   0 &  0 &  0 &    0\\
 0 &  -1 &   0 &  1 &  1 &    0\\
 0 &   0 &  -1 &  1 &  1 &  1/2\\
 0 &   0 &   0 &  1 &  0 &  1/2\\
 0 &   0 &   0 &  0 &  1 &  1/2\\
 0 &   0 &   0 &  0 &  0 &    1\\
\end{bmatrix},
B = 
\begin{bmatrix}
 -1 &  0 &   0 &   0 &   0 & 2/3\\
  0 &  1 &   0 &  -1 &  -1 &   0\\
  0 &  1 &  -1 &   0 &   0 & 1/2\\
  0 &  1 &   0 &   0 &  -1 & 1/2\\
  0 &  1 &   0 &  -1 &   0 &   0\\
  0 &  0 &   0 &   0 &   0 &   1\\
\end{bmatrix}.
\]
Using character theory we get that a faithful representation $\varphi \colon G \to \SO(5,\Z)$, $\R$-equivalent to $id_G$, may be defined by
\[
r(A) \mapsto
\begin{bmatrix}
 1 &   0 &  0 &   0 &  0\\
 0 &  -1 &  0 &   0 &  0\\
 0 &   0 &  1 &   0 &  0\\
 0 &   0 &  0 &  -1 &  0\\
 0 &   0 &  0 &   0 &  1\\
\end{bmatrix},
r(B) \mapsto
\begin{bmatrix}
 1 &   0 &   0 &   0 &  0\\
 0 &  -1 &   0 &   0 &  0\\
 0 &   0 &  -1 &   0 &  0\\
 0 &   0 &   0 &   0 &  1\\
 0 &   0 &   0 &  -1 &  0\\
\end{bmatrix}.
\]

We get that
\[
\lambda_5(\pm e_2e_4) = \varphi r(A), ~\ \lambda_5\left(\pm\frac{e_2e_3(1+e_5e_4)}{\sqrt{2}}\right) = \varphi r(B).
\]

For a map $\epsilon \colon \G \to \Spin(5)$ to be a homomorphism, we have to have the following relations:
\begin{enumerate}
\item The action of $G$ on $\Z^5$:
\begin{equation}
\label{eq:ex:rels1}
\left\{
\begin{aligned}
%\varepsilon(a_1)^2=1,\\
\varepsilon(a_2)\varepsilon(a_3) = 1,\\
\varepsilon(a_2)\varepsilon(a_4)\varepsilon(a_5) = 1.
\end{aligned}
\right.
\end{equation}
\item The relations from $G$. We have
\[
A^2 = 
\begin{bmatrix}
 1 &  0 &  0 &  0 &  0 &  0\\
 0 &  1 &  0 &  0 &  0 &  1\\
 0 &  0 &  1 &  0 &  0 &  1\\
 0 &  0 &  0 &  1 &  0 &  1\\
 0 &  0 &  0 &  0 &  1 &  1\\
 0 &  0 &  0 &  0 &  0 &  1\\
\end{bmatrix},
B^4 = 
\begin{bmatrix}
 1 &  0 &  0 &  0 &  0 &   0\\
 0 &  1 &  0 &  0 &  0 &   0\\
 0 &  0 &  1 &  0 &  0 &   0\\
 0 &  0 &  0 &  1 &  0 &   1\\
 0 &  0 &  0 &  0 &  1 &  -1\\
 0 &  0 &  0 &  0 &  0 &   1\\
\end{bmatrix}
\]
and
\[
(AB)^2 = 
\begin{bmatrix}
 1 &  0 &  0 &  0 &  0 &  0\\
 0 &  1 &  0 &  0 &  0 &  1\\
 0 &  0 &  1 &  0 &  0 &  0\\
 0 &  0 &  0 &  1 &  0 &  1\\
 0 &  0 &  0 &  0 &  1 &  0\\
 0 &  0 &  0 &  0 &  0 &  1\\
\end{bmatrix}.
\]
We get that the following relations should be satisfied in $\Spin(5)$
\begin{equation}
\label{eq:ex:rels2}
\left\{
\begin{aligned}
\varepsilon(A)^2\varepsilon(a_2)\varepsilon(a_3)\varepsilon(a_4)\varepsilon(a_5) = 1\\
\varepsilon(B)^4\varepsilon(a_4)\varepsilon(a_5) = 1\\
(\varepsilon(A)\varepsilon(B))^2\varepsilon(a_2)\varepsilon(a_4) = 1
\end{aligned}
\right.
\end{equation}
\end{enumerate}

From \eqref{eq:ex:rels1} and \eqref{eq:ex:rels2} we get the following conditions on values of $\varepsilon$:
\[
\left\{
\begin{aligned}
%\varepsilon(a_1)^2 & = 1\\
\varepsilon(a_2) & = \varepsilon(a_3) = \varepsilon(a_4)\varepsilon(a_5) = \varepsilon(A)^2 = \varepsilon(B)^4\\
\varepsilon(a_5) & = (\varepsilon(A)\varepsilon(B))^2
%\varepsilon(a_5) & = (\varepsilon(A)\varepsilon(B))^2\\
%\varepsilon(a_4) & = \varepsilon(a_5)\varepsilon(A)^2\\
%\varepsilon(a_1) & = \varepsilon(a_2) = \varepsilon(a_3) = \varepsilon(a_4)\varepsilon(a_5)
\end{aligned}
\right.
\]
Note that both values of $\varepsilon(a_1)$ are allowed. Recall that
\[
\varepsilon(A) = \pm e_2e_4 \text{ and } \varepsilon(B) = \pm\frac{e_2e_3(1+e_5e_4)}{\sqrt{2}}.
\]
Since for any of the above values we have $\varepsilon(A)^2 = \varepsilon(B)^4 = (\varepsilon(A)\varepsilon(B))^2 = -1$, hence we get $8$ spin structures on $\R^5/\G$ and there exists a spin structure on $\R^5/\G'$. Moreover, since $H^1(\R^5/\G',\Z_2) = \Z_2^2$, we get exactly four spin structures on the former manifold.

\section{Some statistics}

\begin{table} %[ht]
\begin{center}
{
\footnotesize
\renewcommand{\arraystretch}{1.2}
\subfloat{
\begin{tabular}{l|l|l|r}
 $\G'$ & $G'$ & $r^{-1}(G) $ & \#S \\ \hline
 min.58.1.1.0   & $1$            & min.58.1.1.0   & 32\\
 min.59.1.1.1   & $ C_2 $ & min.59.1.1.1   & 32\\
 min.62.1.1.1   & $ C_2 $ & min.62.1.1.1   & 32\\
 min.62.1.2.1   & $ C_2 $ & min.62.1.2.1   & 16\\
 min.62.1.3.1   & $ C_2 $ & min.62.1.3.1   &  8\\
 min.65.1.1.7   & $(C_2)^2$ & min.65.1.1.7   & 16\\
 min.66.1.1.11  & $(C_2)^2$ & min.66.1.1.11  & 16\\
 min.66.1.3.11  & $(C_2)^2$ & min.66.1.3.11  &  8\\
 min.70.1.1.20  & $(C_2)^2$ & min.70.1.1.20  & 32\\
 min.70.1.1.22  & $(C_2)^2$ & min.70.1.1.22  & 16\\
 min.70.1.1.28  & $(C_2)^2$ & min.70.1.1.28  & 16\\
 min.70.1.1.30  & $(C_2)^2$ & min.70.1.1.30  & 16\\
 min.70.1.14.1  & $(C_2)^2$ & min.70.1.14.1  &  4\\
 min.70.1.15.19 & $(C_2)^2$ & min.70.1.15.19 & 16\\
 min.70.1.15.5  & $(C_2)^2$ & min.70.1.15.5  &  8\\
 min.70.1.1.76  & $(C_2)^2$ & min.70.1.1.76  & 16\\
 min.70.1.1.77  & $(C_2)^2$ & min.70.1.1.77  & 16\\
 min.70.1.1.94  & $(C_2)^2$ & min.70.1.1.94  & 16\\
 min.70.1.2.25  & $(C_2)^2$ & min.70.1.2.25  &  8\\
 min.70.1.2.9   & $(C_2)^2$ & min.70.1.2.9   & 16\\
 min.70.1.3.11  & $(C_2)^2$ & min.70.1.3.11  &  8\\
 min.70.1.3.7   & $(C_2)^2$ & min.70.1.3.7   & 16\\
 min.70.1.4.10  & $(C_2)^2$ & min.70.1.4.10  &  8\\
 min.70.1.4.11  & $(C_2)^2$ & min.70.1.4.11  &  8\\
 min.70.1.4.7   & $(C_2)^2$ & min.70.1.4.7   & 16\\
 min.70.1.4.9   & $(C_2)^2$ & min.70.1.4.9   &  8\\
 min.70.1.6.3   & $(C_2)^2$ & min.70.1.6.3   &  8\\
 min.70.1.7.13  & $(C_2)^2$ & min.70.1.7.13  &  8\\
 min.70.1.7.15  & $(C_2)^2$ & min.70.1.7.15  &  8\\
 min.71.1.1.362 & $(C_2)^3$ & min.71.1.1.362 & 16\\
 min.71.1.1.371 & $(C_2)^3$ & min.71.1.1.371 &  8\\
 min.71.1.1.373 & $(C_2)^3$ & min.71.1.1.373 & 16\\
 min.71.1.1.375 & $(C_2)^3$ & min.71.1.1.375 &  8\\
 min.71.1.1.378 & $(C_2)^3$ & min.71.1.1.378 &  8\\
 min.71.1.1.382 & $(C_2)^3$ & min.71.1.1.382 &  8\\
 min.71.1.25.95 & $(C_2)^3$ & min.71.1.25.95 & 16\\
 min.75.1.1.1   & $C_4$ & min.75.1.1.1   &  8\\
 min.79.1.1.1   & $C_4$ & min.79.1.1.1   & 16\\
 min.79.1.2.2   & $C_4$ & min.79.1.2.2   &  8\\
 min.81.1.1.1   & $C_4$ & min.81.1.1.1   & 16\\
 min.81.1.3.1   & $C_4$ & min.81.1.3.1   &  8\\
 min.81.1.6.1   & $C_4$ & min.81.1.6.1   &  8\\
 min.85.1.1.41  & $D_8$ & min.85.1.1.41  & 16\\
 min.85.1.1.42  & $D_8$ & min.85.1.1.42  &  8\\
\end{tabular}
}
\subfloat{
\begin{tabular}{l|l|l|r}
 $\G'$ & $G'$ & $r^{-1}(G) $ & \#S \\ \hline
 min.85.1.1.44  & $D_8$ & min.85.1.1.44  &  8\\
 min.85.1.1.45  & $D_8$ & min.85.1.1.45  &  8\\
 min.85.1.1.46  & $D_8$ & min.85.1.1.46  &  8\\
 min.85.1.3.19  & $D_8$ & min.85.1.3.19  &  4\\
 min.85.1.3.22  & $D_8$ & min.85.1.3.22  &  8\\
 min.86.1.13.5  & $D_8$ & min.86.1.13.5  &  4\\
 min.86.1.13.6  & $D_8$ & min.86.1.13.6  &  8\\
 min.86.1.13.7  & $D_8$ & min.86.1.13.7  &  4\\
 min.90.1.10.3  & $ C_4 \times C_2$ & min.90.1.10.3  &  4\\
 min.98.1.3.12  & $ C_4 \times C_2$ & min.98.1.3.12  &  4\\
 min.101.1.1.1  & $ C_3$ & min.58.1.1.0   &  2\\
 min.104.1.1.1  & $ C_3$ & min.58.1.1.0   &  8\\
 min.104.1.2.1  & $ C_3$ & min.58.1.1.0   &  8\\
 min.106.1.1.1  & $ C_6$ & min.62.1.1.1   &  8\\
 min.107.1.1.2  & $ S_3 $ & min.62.1.2.1   &  8\\
 min.107.1.2.1  & $ S_3 $ & min.62.1.3.1   &  4\\
 min.107.2.1.2  & $ S_3 $ & min.62.1.2.1   &  8\\
 min.107.2.2.1  & $ S_3 $ & min.62.1.3.1   &  4\\
 min.107.2.3.2  & $ S_3 $ & min.62.1.2.1   &  8\\
 min.107.2.4.1  & $ S_3 $ & min.62.1.3.1   &  4\\
 min.110.1.1.1  & $ C_6 $ & min.59.1.1.1   &  8\\
 min.110.1.3.1  & $ C_6 $ & min.59.1.1.1   &  8\\
 min.123.1.1.1  & $ C_{12}       $ & min.79.1.1.1   &  4\\
 min.124.1.1.1  & $ C_{12}       $ & min.81.1.1.1   &  4\\
 min.129.1.1.1  & $ C_6 $ & min.62.1.1.1   &  2\\
 min.129.1.2.1  & $ C_6 $ & min.62.1.3.1   &  2\\
 min.130.1.1.12 & $ (C_3)^2$ & min.58.1.1.0   &  2\\
 min.130.1.1.37 & $ (C_3)^2$ & min.58.1.1.0   &  2\\
 min.130.1.3.10 & $ (C_3)^2$ & min.58.1.1.0   &  2\\
 min.131.1.2.3  & $ A_4        $ & min.70.1.15.19 &  4\\
 min.131.2.1.3  & $ A_4        $ & min.70.1.1.76  &  4\\
 min.132.1.2.3  & $ C_2 \times A_4   $ & min.71.1.25.95 &  4\\
 min.132.2.1.6  & $ C_2 \times A_4   $ & min.71.1.1.373 &  4\\
 min.134.1.2.2  & $ S4        $ & min.86.1.13.6  &  4\\
 min.144.1.1.1  & $C_8$  & min.144.1.1.1 & 4\\
 min.154.1.1.1  & $ C_{12}       $ & min.75.1.1.1   &  2\\
 min.164.1.1.1  & $ C_5        $ & min.58.1.1.0   &  2\\
 group.240.2.1.11 & $ D_8        $ & group.240.2.1.11 &  8\\
 group.326.1.1.1  & $ C_6        $ & min.59.1.1.1     &  2\\
 group.341.1.1.1  & $ C_6        $ & min.62.1.1.1     &  8\\
 group.361.1.1.21 & $ D_{12}       $ & min.70.1.3.7     &  8\\
 group.361.1.1.22 & $ D_{12}       $ & min.70.1.3.11    &  4\\
 group.541.1.1.10 & $ C_6 \times C_3   $ & min.62.1.1.1     &  2\\
 group.994.1.1.1  & $ C_{10}       $ & min.59.1.1.1     &  2\\
\end{tabular}
}
}
\end{center}
\caption{Spin structures in dimension 5. $\G'$ -- the name of the Bieberbach group, $G'$ -- isomorphism type of the holonomy group of $\G'$, $r^{-1}(G)$ -- the preimage of the Sylow 2-subgroup of $G$, \#S -- the number of spin structures on $\R^5/\G'$}
\label{tab:spin5}
\end{table}

Recall that CARAT represents any $n$-dimensional Bieberbach group $\G$ as subgroup of $\GL(n,\Z) \ltimes \Q^n$. In this representation the maximal normal abelian subgroup equals $\Z^n$, the holonomy group $G \cong \G/\Z^n$ is a finite subgroup of $\GL(n,\Z)$ and the integral holonomy representation is just the inclusion map to $\GL(n,\Z)$. By the $\Z$-class and the $\Q$-class of a finite subgroup of $\GL(n,\Z)$ we will denote the conjugacy class of the group in $\GL(n,\Z)$ and $\GL(n,\Q)$ respectively.

A necessary condition for a Bieberbach group to have a spin structure is to be orientable. This property is fully determined by the $\Q$-class of the holonomy group. Note that the existence of spin structures is -- in contrast to orientation -- determined only by the isomorphism class of a Bieberbach group (see Remark \ref{rem:dependency} below). Table \ref{tab:qcl} shows the number of $\Q$-classes, the number of $\Q$-classes which determine orientable flat manifolds and the number of $\Q$-classes for which there exists a flat manifold with a spin structure in dimensions 5 and 6.

\begin{table} %[h]

\centering
\begin{tabular}{r|r|r|r}
Dim & \#$\Q$C & \#O$\Q$C & \#S$\Q$C\\ \hline
5 & 95 & 41 & {35}\\
6 &397 &106 & {92}\\
\end{tabular}

\caption{Number of all $\Q$-classes ($\Q$C), orientable $\Q$-classes (O$\Q$C) and the ones for which there exist spin manifolds (S$\Q$C) in dimensions 5 and 6}
\label{tab:qcl}
\end{table}

Table \ref{tab:bg} shows the number of all flat manifolds, orientable flat manifolds and flat manifolds which admit a spin structure in dimensions 5 and 6. Because of their number, in Table \ref{tab:spin5} we list all Bieberbach groups with spin structures of dimension 5. The data for both dimensions 5 and 6 can downloaded from the WWW page \cite{LPdata}.

\begin{table} %[ht]

\centering
\begin{tabular}{r|r|r|r}
Dim & \#FM & \#OFM & \#SFM\\ \hline
5 & 1060 & 174 & {88}\\
6 & 38746& 3314& {760}\\
\end{tabular}

\caption{Number of all, oriented and spin flat manifolds in dimensions 5 and 6}
\label{tab:bg}
\end{table}

\begin{rem}
\label{rem:dependency}
From the paper \cite{PS10} the following facts hold for flat manifolds in dimension 4:
\begin{enumerate}
\item The existence of a spin structure does not depend on the $\Q$-class of the integral holonomy representation of an orientable flat manifold.
\item The existence of a spin structure is determined by the $\Z$-class of the integral holonomy representation of an orientable flat manifold.
\end{enumerate}
By \cite[Example 3.3]{HiSz08} and \cite[Theorem 3.2]{MP04} we know that the former fact does not hold in dimension 6.
The calculations give 5-dimensional examples -- for each of the following $\Z$-classes of finite subgroups of $\GL(5,\Z)$ there exist Bieberbach groups with holonomy group in the class, with and without spin structures:

\begin{center}
\begin{tabular}{llllll}
min.66.1.1,& min.66.1.3,& min.70.1.1,& min.70.1.15,& min.70.1.2,& min.70.1.3,\\ min.70.1.7,& min.71.1.1,& min.71.1.25,& min.85.1.3, & \multicolumn{2}{l}{group.361.1.1.}
\end{tabular}
\end{center}
Those $\Z$-classes belong to the following 5 $\Q$-classes:
\begin{center}
\begin{tabular}{lllll}
min.66,& min.70, & min.71,& min.85, & group.361.
\end{tabular}
\end{center}

Moreover there are {100} $\Z$-classes of finite subgroups of $\GL(6,\Z)$, collected in {37} $\Q$-classes, for which we can find examples of Bieberbach groups with and without spin structures.

\end{rem}

\section*{Acknowledgments}

The computations were performed with usage of Maxima \cite{maxima}, Carat \cite{CARAT03} and GAP \cite{GAP14}, in particular GAP package HAP \cite{HAP13}.

This article was supported by the National Science Center Poland grant no. 2013/09/B/ST1/04125.

%\section*{References}
\bibliographystyle{plain}
\bibliography{bibl}

\newpage

\appendix

\section{Corrigendum to ``Spin structures on flat manifolds''}

\noindent
We underline the changes which has been made in order to present corrected results of the paper.

\subsection{Introduction}

Due to a computational oversight stemming from an error in code that is no longer accessible, we erroneously stated in the original article \cite{LP15} that all Bieberbach groups within the following $\mathbb{Q}$-classes are fundamental groups of flat manifolds without spin structures:
\begin{center}
\begin{tabular}{llllllll}
min.141, & min.207, & min.264, & min.265, & min.268, & min.270, & min.468, & group.1264.
\end{tabular}
\end{center}
However, within this set, one five-dimensional group and 43 six-dimensional groups do in fact admit spin structures.

Detailed instructions for reproducing our results are available in the repository \cite{Lu25GitHub}.

\subsection{Corrections}

The calculation errors affected only the results presented in Section 8 of \cite{LP15}. These results are discussed in the following subsections, with the exception of the last sentence on page 290, which should be replaced with the following corrected version:

\vspace{.5em}
``Moreover there are \rl{100} $\Z$-classes of finite subgroups of $\GL(6,\Z)$, collected in \rl{37} $\Q$-classes, for which we can find examples of Bieberbach groups with and without spin structures.''
\vspace{.5em}

Note that we use boldface to indicate corrections or newly added material.

\subsection{Table 1}

In Table 1 there should be one row added:
\begin{center}
\begin{tabular}{l|l|l|r}
 $\G'$ & $G'$ & $r^{-1}(G) $ & \#S \\ \hline
\rl{min.144.1.1.1} & $\mathbf{C_8}$ & \rl{min.144.1.1.1} & \rl{4}\\
\end{tabular}
\end{center}

\subsection{Table 2}

Table 2 should have the following form:

\begin{center}
\begin{tabular}{r|r|r|r}
Dim & \#$\Q$C & \#O$\Q$C & \#S$\Q$C\\ \hline
5 & 95 & 41 & \rl{35}\\
6 &397 &106 & \rl{92}\\
\end{tabular}
\end{center}

\subsection{Table 3}

Table 2 should have the following form:

\begin{center}
\begin{tabular}{r|r|r|r}
Dim & \#FM & \#OFM & \#SFM\\ \hline
5 & 1060 & 174 & \rl{88}\\
6 & 38746& 3314&\rl{760}\\
\end{tabular}
\end{center}

\section*{Acknowledgments}

The authors would like to thank Miguel Montero for pointing out the omission of the five-dimensional group. His inquiry prompted a recalculation of spin structures on low-dimensional flat manifolds, which led to the discovery of the aforementioned errors.

\end{document}